\documentclass[12pt]{article}
\usepackage{amsmath}
\usepackage{amssymb,amsfonts,amsmath,amsthm}
\usepackage{epsfig}
\newtheorem{theorem}{Theorem}[section]
\newtheorem{thm}[theorem]{Theorem}
\newtheorem{prop}[theorem]{Proposition}
\newtheorem{lem}[theorem]{Lemma}

\newtheorem{cor}[theorem]{Corollary}

\makeatletter \@addtoreset{equation}{section}

\def\CT{\mathop{\mathrm{CT}}}

\begin{document}

\title{Two Coefficients of the Dyson Product}
\author{
  {\vspace{0.2cm}}
  {Lun Lv$^1$, Guoce Xin$^2$, and Yue Zhou$^3$}\\
  {\small $^{1,2,3}$Center for Combinatorics, LPMC--TJKLC}\\
  {\small Nankai University,}\\
  {\small Tianjin 300071, P.R. China}\\
  { \small $^1$lvlun@cfc.nankai.edu.cn\ \ \ $^2$gxin@nankai.edu.cn\ \ \ $^3$nkzhouyue@gmail.com}\\
 }

\date{November  21, 2007}
\maketitle

\begin{abstract}
In this paper, the closed-form expressions for the coefficients of
$\frac{x_r^2}{x_s^2}$ and $\frac{x_r^2}{x_sx_t}$ in the Dyson
product are found by applying an extension of Good's idea. As
consequences, we find several interesting Dyson style constant term
identities.
\end{abstract}



\section{Introduction}

For nonnegative integers $a_1,a_2,\ldots,a_n$, define
\begin{align*}
\qquad \qquad \qquad \quad\quad \quad  \quad  D_n(\mathbf{x},\mathbf{a}):= \prod_{1\le i\ne j \le n}
\left(1-\frac{x_i}{x_j}\right)^{\!\!a_i},
\qquad  \qquad\,\quad(\mbox{Dyson product})
\end{align*}
where $\mathbf{x}:=(x_1,\ldots,x_n)$ and $\mathbf{a}:=(a_1,\ldots,a_n)$.

Dyson \cite{dyson1962} conjectured the following constant term identity in 1962.
\begin{thm}[Dyson's Conjecture]\label{dyson}
\begin{align*}
\CT_{\mathbf{x}} D_n(\mathbf{x},\mathbf{a}) =
\frac{(a_1+a_2+\cdots+a_n)!}{ a_1!\,a_2!\,\cdots\, a_n!}.
\end{align*}
where $\CT_{\mathbf{x}}f(\mathbf{x})$ means to take the
constant term in the $x$'s of the series $f(\mathbf{x})$.
\end{thm}

Dyson's conjecture was first proved independently by
Gunson~\cite{gunson} and  by Wilson \cite{wilson}. Later an elegant
recursive proof was published by Good \cite{good1}. Andrews
\cite{andrews1975} conjectured the $q$-analog of the Dyson
conjecture which was first proved, combinatorially, by Zeilberger
and Bressoud \cite{zeil-bres1985} in 1985. Recently, Gessel and Xin
\cite{gess-xin2006} gave a very different proof by using properties
of formal Laurent series and of polynomials.

Good's idea has been extended by several authors. The current
interest is to evaluate the coefficients of monomials $M$ of degree
$0$ in the Dyson product, where $M:=\prod_{i=0}^nx_{i}^{b_i}$ with
$\sum_{i=0}^nb_i=0$. Kadell \cite{kadell1998} outlined the use of
Good's idea for $M$ to be
$\frac{x_1}{x_n},\frac{x_1x_2}{x_{n-1}x_n}$ and
$\frac{x_1x_2}{x_n^2}$. Along this line, Zeilberger and Sills
\cite{sills-zeilberger} presented a case study in experimental yet
rigorous mathematics by describing an algorithm that automatically
conjectures and proves closed-form. Using this algorithm, Sills
\cite{sills2006} guessed and proved closed-form expressions for $M$
to be $\frac{x_s}{x_r}$, $\frac{x_sx_t}{x_r^2}$ and
$\frac{x_tx_u}{x_rx_s}$. These results and their $q$-analogies were
recently generalized for $M$ with a square free numerator by Lv, Xin
and Zhou \cite{lv_xin_zhou} by extending Gessel-Xin's Laurent series
method \cite{gess-xin2006} for proving the $q$-Dyson Theorem.

The cases for $M$ having a square in the numerator are much more
complicated. By extending Good's idea, we obtain closed forms for
the simplest cases $M=\frac{x_r^2}{x_s^2}$ and
$M=\frac{x_r^2}{x_sx_t}$. In doing so, we guess these two formulas
simultaneously, written as a sum instead of a single product. Our
main results are stated as follows.

\begin{thm}\label{2dy1}
Let $r$ and $s$ be distinct integers with $1\leq r, s\leq n$. Then
{\small\begin{align}\label{main1} \CT_{\mathbf{x}}
\frac{x_s^2}{x_r^2}D_n(\mathbf{x},\mathbf{a})
=\frac{a_r}{(1+a^{(r)})(2+a^{(r)})}\bigg[(a_r-1)-\sum_{\substack{i=1\\i\neq
r,s}}^{n} \frac{a_i(1+a)}{(1+a^{(r)}-a_i)}\bigg]C_n(\mathbf{a}),
\end{align}}
where $a:=a_1+a_2+\cdots+a_n$, $a^{(j)}:=a-a_j$ and
$C_n(\mathbf{a}):=\frac{(a_1+a_2+\cdots+a_n)!}{ a_1!\,a_2!\,\cdots\, a_n!}$.
\end{thm}

\begin{thm}\label{2dy2}
Let $r,s$ and $t$ be distinct integers with $1\leq r, s,t\leq n$.
Then
{\small\begin{align}\label{main2} \CT_{\mathbf{x}}
\frac{x_sx_t}{x_r^2}D_n(\mathbf{x},\mathbf{a})
=\frac{a_r}{(1+a^{(r)})(2+a^{(r)})}\bigg[(a+a_r)-\sum_{\substack{i=1\\i\neq
r,s,t}}^{n} \frac{a_i(1+a)}{(1+a^{(r)}-a_i)}\bigg]C_n(\mathbf{a}),
\end{align}}
where $a, a^{(r)}$ and $C_n(\mathbf{a})$ are defined as Theorem $\mathrm{\ref{2dy1}}$.
\end{thm}

The proofs will be given in Section 2. In Section 3, we construct
several interesting Dyson style constant term identities.

\section{Proof of Theorem \ref{2dy1} and  Theorem \ref{2dy2}}
Good's proof \cite{good1} of the Dyson Conjecture uses the
recurrence
\begin{align*}
D_n(\mathbf{x},\mathbf{a})=\sum_{k=1}^n
D_n(\mathbf{x},\mathbf{a}-\mathbf{e}_k),
\end{align*}
where $\mathbf{e}_k:=(0,\ldots,0,1,0,\ldots,0)$ is the $k$th unit
coordinate $n$-vector. It follows that the following recurrence
holds for any monomial $M$ of degree $0$.
\begin{align*}
\CT_{\mathbf{x}} \frac{1}{M}D_n(\mathbf{x},\mathbf{a})=\sum_{k=1}^n \CT_{\mathbf{x}}
\frac{1}{M} D_n(\mathbf{x},\mathbf{a}-\mathbf{e}_k).
\end{align*}
Thus if we can guess a formula, then we can prove it by checking the
initial condition, the recurrence and the boundary conditions. This
is the so called Good-style proof.

In our case, denote by $F_L(r,s,\mathbf{a})$ (resp. $
G_L(r,s,t,\mathbf{a})$) the left-hand side of \eqref{main1} (resp.
\eqref{main2}), and by $F_R(r,s,\mathbf{a})$ (resp.
$G_R(r,s,t,\mathbf{a})$) the right-hand side of \eqref{main1} (resp.
\eqref{main2}). Without loss of generality, we may assume $r=1,s=2$
and $t=3$ in Theorems \ref{2dy1} and \ref{2dy2}, i.e., we need to prove that
\begin{align*}
F_L(\mathbf{a})=F_R(\mathbf{a}),\quad\quad G_L(\mathbf{a})=G_R(\mathbf{a}),
\end{align*}
where  $F_L(\mathbf{a}):=F(1,2,\mathbf{a})$ and we use similar
notations for $F_R(\mathbf{a}),G_L(\mathbf{a})$ and $G_R(\mathbf{a})$.

\subsection{Initial Condition}

We can easily verify that
\begin{align*}
&F_L(\mathbf{0})=F_R(\mathbf{0})=0,\\
&G_L(\mathbf{0})=G_R(\mathbf{0})=0.
\end{align*}

\subsection{Recurrence}

We need to show that $F_R(\mathbf{a})$ and $G_R(\mathbf{a})$ satisfy
the recurrences
\begin{align}
F_R(\mathbf{a})=\sum_{k=1}^n F_R(\mathbf{a}-\mathbf{e}_k),
\label{recurr1}\\
G_R(\mathbf{a})=\sum_{k=1}^n G_R(\mathbf{a}-\mathbf{e}_k).
\label{recurr2}
\end{align}
In order to do so, we define
\begin{align*}
&H_1(\mathbf{a}):=\frac{a_1(a_1-1)}{(1+a^{(1)})(2+a^{(1)})}C_n(\mathbf{a}),\\
&H_2(\mathbf{a}):=\frac{a_1(a+a_1)}{(1+a^{(1)})(2+a^{(1)})}C_n(\mathbf{a}),\\
&H_i(\mathbf{a}):=\frac{a_1a_i(1+a)}{(1+a^{(1)})(2+a^{(1)})(1+a^{(1)}-a_i)}C_n(\mathbf{a}),\
\ i=3,4,\ldots,n.
\end{align*}
Then
$F_R(\mathbf{a})=H_1(\mathbf{a})+\sum_{i=3}^nH_i(\mathbf{a})$
and
$G_R(\mathbf{a})=H_2(\mathbf{a})+\sum_{i=4}^nH_i(\mathbf{a})$.
Therefore to show \eqref{recurr1} and \eqref{recurr2}, it suffices
to show the following:

\begin{lem}
For each $i=1,2,\ldots,n$, we have the recurrence
$H_i(\mathbf{a})=\sum_{k=1}^n H_i(\mathbf{a}-\mathbf{e}_k)$.
\end{lem}

\begin{proof}
\begin{enumerate}
\item
For $H_1(\mathbf{a})$,
{\small\begin{align*}
\sum_{k=1}^n H_1(\mathbf{a}-\mathbf{e}_k)&=\frac{(a_1-1)(a_1-2)}{(1+a^{(1)})(2+a^{(1)})}C_n(\mathbf{a}-\mathbf{e}_1)
  +\sum_{k=2}^{n}\frac{a_1(a_1-1)}{a^{(1)}(1+a^{(1)})}C_n(\mathbf{a}-\mathbf{e}_k)\\
&=\bigg[\frac{a_1(a_1-1)(a_1-2)}{a(1+a^{(1)})(2+a^{(1)})}
  +\sum_{k=2}^{n}\frac{a_ka_1(a_1-1)}{aa^{(1)}(1+a^{(1)})}\bigg]C_n(\mathbf{a})\\
&=\bigg[\frac{a_1(a_1-1)(a_1-2)}{a(1+a^{(1)})(2+a^{(1)})}
  +\frac{a_1(a_1-1)}{a(1+a^{(1)})}\bigg]C_n(\mathbf{a})\\
&=\frac{a_1(a_1-1)}{(1+a^{(1)})(2+a^{(1)})}C_n(\mathbf{a})=H_1(\mathbf{a}).
\end{align*}}
\item
For $H_2(\mathbf{a})$,
{\small\begin{align*}
\sum_{k=1}^n H_2(\mathbf{a}-\mathbf{e}_k)&=\frac{(a_1-1)(a+a_1-2)}{(1+a^{(1)})(2+a^{(1)})}C_n(\mathbf{a}-\mathbf{e}_1)
  +\sum_{k=2}^{n}\frac{a_1(a+a_1-1)}{a^{(1)}(1+a^{(1)})}C_n(\mathbf{a}-\mathbf{e}_k)\\
&=\bigg[\frac{a_1(a_1-1)(a+a_1-2)}{a(1+a^{(1)})(2+a^{(1)})}
  +\sum_{k=2}^{n}\frac{a_ka_1(a+a_1-1)}{aa^{(1)}(1+a^{(1)})}\bigg]C_n(\mathbf{a})\\
&=\bigg[\frac{a_1(a_1-1)(a+a_1-2)}{a(1+a^{(1)})(2+a^{(1)})}
  +\frac{a_1(a+a_1-1)}{a(1+a^{(1)})}\bigg]C_n(\mathbf{a})\\
&=\frac{a_1(a+a_1)}{(1+a^{(1)})(2+a^{(1)})}C_n(\mathbf{a})=H_2(\mathbf{a}).
\end{align*}}
\item
For $H_i(\mathbf{a})$ with $i=3,\ldots,n$, without loss of generality, we may assume $i=3$.
{\small\begin{align*}
&\sum_{k=1}^n H_3(\mathbf{a}-\mathbf{e}_k)\\
=&\frac{aa_3(a_1-1)}
  {(1+a^{(1)})(2+a^{(1)})(1+a^{(1)}-a_3)}C_n(\mathbf{a}-\mathbf{e}_1)
  +\frac{aa_1a_3}{a^{(1)}(1+a^{(1)})(a^{(1)}-a_3)}C_n(\mathbf{a}-\mathbf{e}_2)\\
&+\frac{aa_1(a_3-1)}{a^{(1)}(1+a^{(1)})(1+a^{(1)}-a_3)}C_n(\mathbf{a}-\mathbf{e}_3)
  +\sum_{k=4}^{n}\frac{aa_1a_3}{a^{(1)}(1+a^{(1)})(a^{(1)}-a_3)}C_n(\mathbf{a}-\mathbf{e}_k)\\
=&\frac{a_1a_3(a_1-1)}
  {(1+a^{(1)})(2+a^{(1)})(1+a^{(1)}-a_3)}C_n(\mathbf{a})
  +\frac{a_1a_2a_3}{a^{(1)}(1+a^{(1)})(a^{(1)}-a_3)}C_n(\mathbf{a})\\
&+\frac{a_1a_3(a_3-1)}{a^{(1)}(1+a^{(1)})(1+a^{(1)}-a_3)}C_n(\mathbf{a})
  +\frac{a_1a_3(a-a_1-a_2-a_3)}{a^{(1)}(1+a^{(1)})(a^{(1)}-a_3)}C_n(\mathbf{a})\\
=&\frac{a_1a_3}{(1+a^{(1)})(1+a^{(1)}-a_3)}\bigg[\frac{a_1-1}{2+a^{(1)}}+\frac{a_3-1}{a^{(1)}}
   +\frac{(a-a_1-a_3)(1+a^{(1)}-a_3)}{a^{(1)}(a^{(1)}-a_3)}\bigg]C_n(\mathbf{a})\\
=&\frac{a_1a_3(1+a)}{(1+a^{(1)})(2+a^{(1)})(1+a^{(1)}-a_3)}C_n(\mathbf{a})=H_3(\mathbf{a}).
\end{align*}}
\end{enumerate}
This completes the proof.
\end{proof}

\subsection{Boundary Conditions}

Now we consider the boundary conditions. For any $k$ with $1\leq
k\leq n$,
{\small\begin{align*}
D_n\big(\mathbf{x},(a_1,\ldots,a_{k-1},0,a_{k+1},\ldots,a_n)\big)
=D_{n-1}(\mathbf{x}^{\langle k \rangle},\mathbf{a}^{\langle k
\rangle}) \times \prod_{\substack{i=1\\i \neq
k}}^{n}\left(1-\frac{x_i}{x_k}\right)^{\!\!a_i},
\end{align*}}
where $\mathbf{x}^{\langle k \rangle}:=(x_1,\ldots,x_{k-1},x_{k+1},\ldots,x_n).$
Thus we have
\begin{align}
\CT_{\mathbf{x}}\ \frac{x_2^2}{x_1^2}D_n\big(\mathbf{x},(a_1,\ldots,a_{k-1},0,a_{k+1},\ldots,a_n)\big)%\nonumber\\
&=\CT_{\mathbf{x}^{\langle k \rangle}}P_k\cdot
D_{n-1}(\mathbf{x}^{\langle k \rangle},\mathbf{a}^{\langle k
\rangle}),\label{coeff1-5}
\end{align}
where $P_k$ is given by
\begin{align*}
P_k:&=\CT_{x_k}\frac{x_2^2}{x_1^2} \prod_{\substack{i=1\\ i \neq
k}}^{n}\Big(1-\frac{x_i}{x_k}\Big)^{\!\!a_i}
\\
&=\left\{%
\begin{array}{ll}
    0, & k=1; \\
    {a_1\choose 2}+a_1\sum_{i=3}^na_i\frac{x_i}{x_1}+\sum_{i=3}^n{a_i\choose 2}\frac{x_i^2}{x_1^2}
    +\sum\limits_{3\leq i<j\leq n}a_ia_j\frac{x_ix_j}{x_1^2}, & k=2; \\
    \frac{x_2^2}{x_1^2}, & \hbox{otherwise.} \\
\end{array}%
\right.
\end{align*}
Taking the constant term in the $x$'s of \eqref{coeff1-5}, we obtain
{\small\begin{align*}
F_L&(a_1,\ldots,a_{k-1},0,a_{k+1},\ldots,a_n)\\
&=\left\{%
\begin{array}{ll}
    0, & k=1; \\
    \CT\limits_{\mathbf{x}^{\langle 2 \rangle}}
    \Big({a_1\choose 2}+a_1\sum\limits_{i=3}^na_i\frac{x_i}{x_1}+\sum\limits_{i=3}^n{a_i\choose 2}\frac{x_i^2}{x_1^2}
    +\sum\limits_{3\leq i<j\leq n}a_ia_j\frac{x_ix_j}{x_1^2}\Big)
    D_{n-1}(\mathbf{x}^{\langle 2 \rangle},\mathbf{a}^{\langle 2 \rangle}), & k=2; \\
    \CT\limits_{\mathbf{x}^{\langle k \rangle}}
    \frac{x_2^2}{x_1^2}D_{n-1}(\mathbf{x}^{\langle k \rangle},\mathbf{a}^{\langle k \rangle}), & \hbox{otherwise.} \\
\end{array}%
\right.
\end{align*}}
By Theorem \ref{dyson} and \cite[Theorem 1.1]{sills2006}, we have
\begin{align*}
 \CT_{\mathbf{x}^{\langle 2 \rangle}}{a_1\choose 2}
D_{n-1}(\mathbf{x}^{\langle 2 \rangle},\mathbf{a}^{\langle 2 \rangle})
=&\frac{a_1(a_1-1)}{2}C_{n-1}(\mathbf{a}^{\langle 2 \rangle}),\\
\CT_{\mathbf{x}^{\langle 2 \rangle}}a_1\sum\limits_{i=3}^na_i\frac{x_i}{x_1}
D_{n-1}(\mathbf{x}^{\langle 2 \rangle},\mathbf{a}^{\langle 2 \rangle})
&=-\frac{a_1^2\ (a^{(1)}-a_2)}{1+a^{(1)}-a_2}
C_{n-1}(\mathbf{a}^{\langle 2 \rangle}).
\end{align*}
So we obtain the following boundary conditions (also recurrences)
{\begin{align*}
F_L&(a_1,\ldots,a_{k-1},0,a_{k+1},\ldots,a_n)\nonumber\\
&=\left\{%
\begin{array}{ll}
    0, & k=1; \\
    \Big(\frac{a_1(a_1-1)}{2}-\frac{a_1^2\ a^{(1)}}{1+a^{(1)}}\Big)
C_{n-1}(\mathbf{a}^{\langle 2 \rangle})\\
\qquad+\sum\limits_{i=3}^n{a_i\choose 2}F_L(1,i,\mathbf{a}^{\langle 2 \rangle})
+\sum\limits_{3\leq i<j\leq n}a_ia_jG_L(1,i,j,\mathbf{a}^{\langle 2 \rangle}), & k=2; \\
    F_L(\mathbf{a}^{\langle k \rangle}), &\hbox{otherwise.} \\
\end{array}%
\right.
\end{align*}}

We need to show that $F_R(a_1,\ldots,a_{k-1},0,a_{k+1},\ldots,a_n)$ satisfies
the same boundary conditions. More precisely, the conditions by
replacing all $F_L$ by $F_R$ and all $G_L$ by $G_R$:
{\begin{align}\label{boundary1}
F_R&(a_1,\ldots,a_{k-1},0,a_{k+1},\ldots,a_n)\nonumber\\
&=\left\{%
\begin{array}{ll}
    0, & k=1; \\
    \Big(\frac{a_1(a_1-1)}{2}-\frac{a_1^2\ a^{(1)}}{1+a^{(1)}}\Big)
C_{n-1}(\mathbf{a}^{\langle 2 \rangle})\\
\qquad+\sum\limits_{i=3}^n{a_i\choose 2}F_R(1,i,\mathbf{a}^{\langle 2 \rangle})
+\sum\limits_{3\leq i<j\leq n}a_ia_jG_R(1,i,j,\mathbf{a}^{\langle 2 \rangle}), & k=2; \\
    F_R(\mathbf{a}^{\langle k \rangle}), &\hbox{otherwise.} \\
\end{array}%
\right.
\end{align}}

Similar computation for $\frac{1}{M}=\frac{x_2x_3}{x_1^2}$ yields the boundary conditions:
{\small\begin{align*}
G_L&(a_1,\ldots,a_{k-1},0,a_{k+1},\ldots,a_n)\nonumber\\
&=\left\{%
\begin{array}{ll}
    0, & k=1; \\
    \frac{a_1^2}{1+a^{(1)}}C_{n-1}(\mathbf{a}^{\langle 2 \rangle})-a_3F_L(1,3,\mathbf{a}^{\langle 2 \rangle})
    - \sum_{i=4}^na_iG_L(1,3,i,\mathbf{a}^{\langle 2 \rangle}), & k=2; \\
     \frac{a_1^2}{1+a^{(1)}}C_{n-1}(\mathbf{a}^{\langle 3 \rangle})-a_2F_L(1,2,\mathbf{a}^{\langle 3 \rangle})
    - \sum_{i=4}^na_iG_L(1,2,i,\mathbf{a}^{\langle 3 \rangle}), & k=3; \\
    G_L(\mathbf{a}^{\langle k \rangle}), & \hbox{otherwise,} \\
\end{array}%
\right.
\end{align*}}
so we need to prove the boundary conditions for
$G_R(a_1,\ldots,a_{k-1},0,a_{k+1},\ldots,a_n)$:
{\small\begin{align}\label{boundary2}
G_R&(a_1,\ldots,a_{k-1},0,a_{k+1},\ldots,a_n)\nonumber\\
&=\left\{%
\begin{array}{ll}
    0, & k=1; \\
    \frac{a_1^2}{1+a^{(1)}}C_{n-1}(\mathbf{a}^{\langle 2 \rangle})-a_3F_R(1,3,\mathbf{a}^{\langle 2 \rangle})
    - \sum_{i=4}^na_iG_R(1,3,i,\mathbf{a}^{\langle 2 \rangle}), & k=2; \\
     \frac{a_1^2}{1+a^{(1)}}C_{n-1}(\mathbf{a}^{\langle 3 \rangle})-a_2F_R(1,2,\mathbf{a}^{\langle 3 \rangle})
    - \sum_{i=4}^na_iG_R(1,2,i,\mathbf{a}^{\langle 3 \rangle}), & k=3; \\
    G_R(\mathbf{a}^{\langle k \rangle}), & \hbox{otherwise.} \\
\end{array}%
\right.
\end{align}}

These are summarized by the following lemma.
\begin{lem}
If $a_k=0$ with $k=1,2,\ldots,n$, then
$F_R(a_1,\ldots,a_{k-1},0,a_{k+1},\ldots,a_n)$ satisfies the
boundary conditions \eqref{boundary1} and
$G_R(a_1,\ldots,a_{k-1},0,a_{k+1},\ldots,a_n)$ satisfies the
boundary conditions \eqref{boundary2}.
\end{lem}
\begin{proof}
We only prove the first part for brevity and similarity.

Since the cases $k=1,3,\ldots,n$ are straightforward, we only prove the case $k=2$. Note that during the proof of this
lemma, we have  $a^{(1)}=a_2+a_3+\cdots+a_n=a_3+\cdots+a_n$ because
$a_2=0$.

Since
{\footnotesize\begin{align}\label{bound-1}
\sum_{i=3}^n&{a_i\choose 2}\sum_{\substack{j=3\\j\neq i}}^n \frac{a_j}{1+a^{(1)}-a_j}
   =\sum_{i=3}^{n}\frac{a_i(a_i-1)}{2}\sum_{j=3}^n\Big(\frac{a_j}
   {1+a^{(1)}-a_j}-\frac{a_i}{1+a^{(1)}-a_i}\Big)\nonumber\\
&=\frac{1}{2}\bigg(\sum_{i=3}^n(a_i^2-a_i)\sum_{j=3}^n\frac{a_j}{1+a^{(1)}-a_j}
   -\sum_{i=3}^n\frac{a_i^3-a_i^2}{1+a^{(1)}-a_i}\bigg)\nonumber\\
&=\frac{1}{2}\bigg(\sum_{j=3}^n\frac{a_j}{1+a^{(1)}-a_j}\sum_{i=3}^na_i^2-a^{(1)}\sum_{j=3}^n\frac{a_j}{1+a^{(1)}-a_j}
   -\sum_{i=3}^n\frac{a_i^3-a_i^2}{1+a^{(1)}-a_i}\bigg),
\end{align}}
we have
{\footnotesize\begin{align}\label{coeff1-3}
\sum\limits_{i=3}^n&{a_i\choose 2}F_R(1,i,\mathbf{a}^{\langle 2 \rangle})\nonumber\\
&=\frac{a_1}{(1+a^{(1)})(2+a^{(1)})}\sum_{i=3}^n{a_i\choose 2}\bigg[(a_1-1)
   -\sum_{\substack{j=3\\j\neq i}}^n \frac{a_j(1+a)}{1+a^{(1)}-a_j}\bigg]
   C_{n-1}(\mathbf{a}^{\langle 2 \rangle})\nonumber\\
&=-\frac{a_1(1+a)}{2(1+a^{(1)})(2+a^{(1)})}
\bigg[\sum_{j=3}^n\frac{a_j}{1+a^{(1)}-a_j}\sum_{i=3}^na_i^2
  -a^{(1)}\sum_{j=3}^n\frac{a_j}{1+a^{(1)}-a_j}\nonumber\\
   &\quad-\sum_{i=3}^n\frac{a_i^3-a_i^2}{1+a^{(1)}-a_i}\bigg]C_{n-1}(\mathbf{a}^{\langle 2 \rangle})
   +\frac{a_1(a_1-1)}{2(1+a^{(1)})(2+a^{(1)})}\bigg(\sum_{i=3}^na_i^2-a^{(1)}\bigg)
   C_{n-1}(\mathbf{a}^{\langle 2 \rangle})\quad\mbox{by \eqref{bound-1}}\nonumber\\
&=-\frac{a_1a^{(1)}(a_1-1)}{2(1+a^{(1)})(2+a^{(1)})}C_{n-1}(\mathbf{a}^{\langle 2 \rangle})
   -\lambda \bigg[(1+a)\sum_{j=3}^n\frac{a_j}{1+a^{(1)}-a_j}\sum_{i=3}^na_i^2\nonumber\\
   &\qquad\qquad\qquad\quad\quad\quad\quad\quad\quad\quad\quad-(1+a)\sum_{j=3}^n
   \frac{a_j^3-a_j^2+a_ja^{(1)}}{1+a^{(1)}-a_j}
   -(a_1-1)\sum_{i=3}^na_i^2\bigg],
\end{align}}
where $\lambda:=\frac{a_1}{2(1+a^{(1)})(2+a^{(1)})}C_{n-1}(\mathbf{a}^{\langle 2 \rangle}).$

Observe that
{\small
\begin{align}\label{bound-2}
\sum\limits_{3\leq i<j\leq n}a_ia_j=\frac{1}{2}\bigg[(a^{(1)})^2-\sum_{k=3}^na_k^2\bigg]
\end{align}}
and
{\footnotesize\begin{align}\label{bound-3}
\sum_{3\leq i<j\leq n}&a_ia_j\sum_{\substack{k=3\\k\neq i,j}}^n\frac{a_k}{1+a^{(1)}-a_k}
    =\sum_{3\leq i<j\leq n}a_ia_j\sum_{k=3}^n
    \Big(\frac{a_k}{1+a^{(1)}-a_k}-\frac{a_i}{1+a^{(1)}-a_i}-\frac{a_j}{1+a^{(1)}-a_j}\Big)\nonumber\\
=&\sum_{3\leq i<j\leq n}a_ia_j\sum_{k=3}^n\frac{a_k}{1+a^{(1)}-a_k}-\sum_{3\leq i<j\leq n}\frac{a_i^2a_j}{1+a^{(1)}-a_i}
    -\sum_{3\leq i<j\leq n}\frac{a_ia_j^2}{1+a^{(1)}-a_j}\nonumber\\
=&\frac{1}{2}\sum_{k=3}^n\frac{a_k}{1+a^{(1)}-a_k}\bigg[\big(a^{(1)}\big)^2-\sum_{i=3}^na_i^2\bigg]
    -\sum_{i=3}^n\sum_{\substack{j=3\\j\neq i}}^n\frac{a_i^2a_j}{1+a^{(1)}-a_i}
    \qquad\qquad\qquad\qquad\mbox{ by \eqref{bound-2}}\nonumber\\
=&\frac{1}{2}\sum_{k=3}^n\frac{a_k}{1+a^{(1)}-a_k}\bigg[\big(a^{(1)}\big)^2-\sum_{i=3}^na_i^2\bigg]
    -\sum_{i=3}^n\frac{a_i^2(a^{(1)}-a_i)}{1+a^{(1)}-a_i}.
\end{align}}
Thus we obtain that
{\footnotesize\begin{align}\label{coeff1-4}
\sum\limits_{3\leq i<j\leq n}&a_ia_jG_R(1,i,j,\mathbf{a}^{\langle 2 \rangle})\nonumber\\
&=\frac{a_1}{(1+a^{(1)})(2+a^{(1)})}\sum\limits_{3\leq i<j\leq n}a_ia_j
   \bigg[(a+a_1)-\sum_{\substack{k=3\\k\neq i,j}}^n\frac{a_k(1+a)}{1+a^{(1)}-a_k}\bigg]
   C_{n-1}(\mathbf{a}^{\langle 2 \rangle})\nonumber\\
&=\frac{a_1(a+a_1)}{(1+a^{(1)})(2+a^{(1)})}\sum\limits_{3\leq i<j\leq n}a_ia_j
   C_{n-1}(\mathbf{a}^{\langle 2 \rangle})\nonumber\\
  &\quad\quad\quad -\frac{a_1(1+a)}{(1+a^{(1)})(2+a^{(1)})}\sum\limits_{3\leq i<j\leq n}a_ia_j
   \sum_{\substack{k=3\\k\neq i,j}}^n\frac{a_k}{1+a^{(1)}-a_k}
   C_{n-1}(\mathbf{a}^{\langle 2 \rangle})\nonumber\\
&=\frac{a_1(a+a_1)}{2(1+a^{(1)})(2+a^{(1)})}\bigg[(a^{(1)})^2-\sum_{k=3}^na_k^2\bigg]
  C_{n-1}(\mathbf{a}^{\langle 2 \rangle})
   +\frac{a_1(1+a)}{2(1+a^{(1)})(2+a^{(1)})}\nonumber\\
&\quad\times\bigg[\sum_{k=3}^n
    \frac{a_k}{1+a^{(1)}-a_k}\Big(\sum_{i=3}^na_i^2-\big(a^{(1)}\big)^2\Big)
    +2\sum_{i=3}^n\frac{a_i^2(a^{(1)}-a_i)}{1+a^{(1)}-a_i}\bigg]
    C_{n-1}(\mathbf{a}^{\langle 2 \rangle})\qquad\mbox{by \eqref{bound-3}}\nonumber\\
&=\frac{a_1(a+a_1)(a^{(1)})^2}{2(1+a^{(1)})(2+a^{(1)})}C_{n-1}(\mathbf{a}^{\langle 2 \rangle})
    +\lambda\bigg[(1+a)\sum_{k=3}^n\frac{a_k}{1+a^{(1)}-a_k}\sum_{i=3}^na_i^2\nonumber\\
    &\quad\quad\quad\quad\quad\quad\quad\quad-(a+a_1)\sum_{k=3}^na_k^2-(1+a)
    \sum_{k=3}^n\frac{a_k(a^{(1)})^2-2a_k^2a^{(1)}+2a_k^3}{1+a^{(1)}-a_k}\bigg].
\end{align}}
Observe that
{\small\begin{align}\label{bound-4}
&(1+a)\sum_{j=3}^n\frac{a_j^3-a_j^2+a_ja^{(1)}}{1+a^{(1)}-a_j}+(a_1-1)\sum_{i=3}^na_i^2\nonumber\\
&\quad\quad\quad\quad\quad\quad\quad-(a+a_1)\sum_{k=3}^na_k^2
-(1+a)\sum_{k=3}^n\frac{a_k(a^{(1)})^2-2a_k^2a^{(1)}+2a_k^3}{1+a^{(1)}-a_k}\nonumber\\
=&(1+a)\sum_{i=3}^n\frac{-a_i^3-a_i^2+a_ia^{(1)}-a_i(a^{(1)})^2+2a_i^2a^{(1)}}{1+a^{(1)}-a_i}
   -(1+a)\sum_{i=3}^na_i^2\nonumber\\
=&(1+a)\sum_{i=3}^n\frac{(1+a^{(1)}-a_i)(a_i^2+2a_i-a_ia^{(1)})-2a_i}{1+a^{(1)}-a_i}
   -(1+a)\sum_{i=3}^na_i^2\nonumber\\
=&(1+a)\sum_{i=3}^n(2a_i-a_ia^{(1)})-(1+a)\sum_{i=3}^n\frac{2a_i}{1+a^{(1)}-a_i}\nonumber\\
=&a^{(1)}(1+a)(2-a^{(1)})-(1+a)\sum_{i=3}^n\frac{2a_i}{1+a^{(1)}-a_i}
\end{align}}
and
{\footnotesize\begin{align}\label{bound-5}
&\frac{a_1a^{(1)}(1+a)(2-a^{(1)})}{2(1+a^{(1)})(2+a^{(1)})} -\frac{a_1a^{(1)}(a_1-1)}{2(1+a^{(1)})(2+a^{(1)})}
+\frac{a_1(a+a_1)(a^{(1)})^2}{2(1+a^{(1)})(2+a^{(1)})}+\frac{a_1(a_1-1)}{2}-\frac{a_1^2a^{(1)}}{1+a^{(1)}}\nonumber\\
&=\frac{a_1(a_1-1)}{(1+a^{(1)})(2+a^{(1)})}.
\end{align}}
Therefore by \eqref{coeff1-3}, \eqref{coeff1-4}, \eqref{bound-4} and \eqref{bound-5},
we have
{\small\begin{align*}
&\bigg[\frac{a_1(a_1-1)}{2}-\frac{a_1^2\ a^{(1)}}{1+a^{(1)}}\bigg]
C_{n-1}(\mathbf{a}^{\langle 2 \rangle})+
\sum\limits_{i=3}^n&{a_i\choose 2}F_R(1,i,\mathbf{a}^{\langle 2 \rangle})
+\sum\limits_{3\leq i<j\leq n}a_ia_jG_R(1,i,j,\mathbf{a}^{\langle 2 \rangle})\\
&=F_R(a_1,0,a_{3},\ldots,a_n).
\end{align*}}
That is to say $F_R(a_1,0,a_{3},\ldots,a_n)$ satisfies boundary conditions \eqref{boundary1}.
\end{proof}

\subsection{The Proof}

Now we can prove Theorems \ref{2dy1} and \ref{2dy2}. Without loss of
generality, we may assume $r=1,s=2$ and $t=3$ in Theorems \ref{2dy1}
and \ref{2dy2}.
\begin{proof}[Proof of Theorems {\rm\ref{2dy1}} and {\rm\ref{2dy2}}]

We prove by induction on $n$ for the two theorems simultaneously.
Clearly, \eqref{main1} and \eqref{main2} hold when $n=2,3$. Assume
that \eqref{main1} and \eqref{main2} hold with $n$ replaced by
$n-1$. Thus for $k=1,2,\ldots,n$ \eqref{main1} and \eqref{main2}
give
\begin{align*}
F_L(r,s,\mathbf{a}^{\langle k \rangle})&=F_R(r,s,\mathbf{a}^{\langle k \rangle}),\\
G_L(r,s,t,\mathbf{a}^{\langle k
\rangle})&=G_R(r,s,t,\mathbf{a}^{\langle k \rangle}).
\end{align*}
That is to say $F_L(a_1,\ldots,a_{k-1},0,a_{k+1},\ldots,a_n)$ and
$F_R(a_1,\ldots,a_{k-1},0,a_{k+1},\ldots,a_n)$ (resp.
$G_L(a_1,\ldots,a_{k-1},0,a_{k+1},\ldots,a_n)$ and
$G_R(a_1,\ldots,a_{k-1},0,a_{k+1},\ldots,a_n)$) satisfy the same
boundary conditions. Additionally $F_L(\mathbf{a})$ and
$F_R(\mathbf{a})$ (resp. $G_L(\mathbf{a})$ and $G_R(\mathbf{a})$)
have the same initial condition and recurrence. It follows that
$F_L(\mathbf{a})=F_R(\mathbf{a})$ (resp.
$G_L(\mathbf{a})=G_R(\mathbf{a})$).
\end{proof}

\section{Several Dyson Style Constant Term Identities}

By linearly combining Theorems \ref{2dy1} and \ref{2dy2}, we obtain
simple formulas.

\begin{prop}\label{pro1}
Let $r,s,t,u$, and $v$ be distinct integers in $\{ 1,2,\dots,n \}$.
Then
\begin{align}
&\CT_{\mathbf{x}}\frac{(x_s-x_t)(x_u-x_v)}{x_r^2}D_n(\mathbf{x},\mathbf{a})=0,\label{propo1}\\
&\CT_{\mathbf{x}}\frac{(x_s-x_u)(x_s-x_v)}{x_r^2}D_n(\mathbf{x},\mathbf{a})
=-\frac{a_r(a+1)}{(2+a^{(r)})(1+a^{(r)}-a_s)}C_n(\mathbf{a}),\label{propo2}\\
&\CT_{\mathbf{x}}\frac{(x_s-x_t)^2}{x_r^2}D_n(\mathbf{x},\mathbf{a})
=-\frac{a_r(a+1)}{2+a^{(r)}}\sum_{i=s,t}\frac{1}{1+a^{(r)}-a_i}C_n(\mathbf{a}).\label{propo3}
\end{align}
\end{prop}
It is worth mentioning that \eqref{propo3} follows from
\eqref{propo1} and \eqref{propo2}, since
$$
(x_s-x_u)(x_s-x_v)+(x_t-x_u)(x_t-x_v)=(x_s-x_t)^2+(x_s-x_u)(x_t-x_v)+(x_s-x_v)(x_t-x_u).
$$

A consequence of Proposition \ref{pro1} is the following:
\begin{cor}
Let I:=$\{i_1,i_2,\ldots,i_{2m}\}$ be a $2m$-element subset of $\{1,2,\ldots,n\}$
and $r\leq n$ is a positive integer
with $r\not\in I$. Then we have
\begin{align*}
\CT_{\mathbf{x}}\frac{\big(\sum_{j=1}^{2m}(-1)^jx_{i_{j}}\big)^2}{x_r^2}&D_n(\mathbf{x},\mathbf{a})
=-\frac{a_r(a+1)}{2+a^{(r)}}\sum_{j\in I}\frac{1}{1+a^{(r)}-a_j}C_n(\mathbf{a}).
\end{align*}
\end{cor}

\begin{proof}
Observe that
\begin{align*}
\bigg(&\sum_{j=1}^{2m}(-1)^jx_{i_{j}}\bigg)^2
=\Big[(x_{i_2}-x_{i_1})+(x_{i_4}-x_{i_3})+\cdots+(x_{i_{2m}}-x_{i_{2m-1}})\Big]^2\\
&=(x_{i_2}-x_{i_1})^2+\cdots+(x_{i_{2m}}-x_{i_{2m-1}})^2
  +\sum_{k=1}^{m}\sum_{\substack{l=1\\l\neq k}}^{m}(x_{i_{2k}}-x_{i_{2k-1}})(x_{i_{2l}}-x_{i_{2l-1}}).
\end{align*}
The corollary then follows by \eqref{propo1} and \eqref{propo3}.
\end{proof}

\medskip
\noindent \textbf{Discussions:} As we have seen in the proof, we
need to guess the formulas of $F_R$ and $G_R$ simultaneously. This
is unlike the coefficients for $M=x_sx_t/x_u^2$ and
$M=x_sx_t/(x_ux_v)$, which have reasonable product formulas and are
equal!

The next cases should be $M$ with $x_r^2x_s$ or $x_r^3$ in the
numerator, both having three cases for the denominator. The
difficulty is: guess three coefficients simultaneously; obtain
enough data.

The study of the $q$-analogies of these formulas will be in a
completely different route and will not be discussed in this paper.

\vspace{.2cm} \noindent{\bf Acknowledgments.}
We would like to acknowledge Zeilberger and Sills' Maple package GoodDyson,
which is available from the web sites
\texttt{http://www.math.rutgers.edu/\\\~{}zeilberg} and
\texttt{http://math.georgiasouthern.edu/\~{}asills}.
This work was
supported by the 973 Project, the PCSIRT project of the Ministry of
Education, the Ministry of Science and Technology and the National
Science Foundation of China.

\end{document}